\documentclass[preprint,12pt, a4paper]{elsarticle}

\usepackage{amssymb}
\usepackage{amsthm}

\usepackage{lineno}
\usepackage{hyperref}
\usepackage{ulem}

\usepackage{float}
\usepackage{tikz}
\usetikzlibrary{positioning, arrows, backgrounds, shapes.multipart, calc}
\usepackage{subcaption}
\usepackage{xurl}

\restylefloat{table}

\journal{arXiv}

\begin{document}

\begin{frontmatter}

%
%

\newpage



\title{Version 2.0 - cashocs: A Computational, Adjoint-Based Shape Optimization and Optimal Control Software}


\author[label1]{Sebastian Blauth\corref{cor1}}
\ead{sebastian.blauth@itwm.fraunhofer.de}
\cortext[cor1]{corresponding author}
\address[label1]{Fraunhofer ITWM, Kaiserslautern, Germany}

\begin{abstract}

In this paper, we present version 2.0 of cashocs. Our software automates the solution of PDE constrained optimization problems for shape optimization and optimal control. Since its inception, many new features and useful tools have been added to cashocs, making it even more flexible and efficient. The most significant additions are a framework for space mapping, the ability to solve topology optimization problems with a level-set approach, the support for parallelism via MPI, and the ability to handle additional (state) constraints. In this software update, we describe the key additions to cashocs, which is now even better-suited for solving complex PDE constrained optimization problems.

\end{abstract}


\end{frontmatter}

{\centering
	This is a post-peer-review, pre-copyedit version of an article published in SoftwareX 24. The final version is available online at \url{https://doi.org/10.1016/j.softx.2023.101577}.
}

\section*{Refers to}
Sebastian Blauth, cashocs: A Computational, Adjoint-Based Shape Optimization and Optimal Control Software, SoftwareX, 13 (2021), p. 100646, \url{https://doi.org/10.1016/j.softx.2020.100646} 

\section*{Metadata}
\label{metadata}

\begin{table}[H]
\begin{tabular}{|l|p{6.5cm}|p{6.5cm}|}

\hline
C1 & Current code version & v2.0.12 \\
\hline
C2 & Permanent link to code/repository used for this code version & \url{https://github.com/sblauth/cashocs/releases/tag/v2.0.12} \\
\hline
C3 & Code Ocean compute capsule & na \\
\hline
C4 & Legal Code License   & GNU GPL v3.0 (or later) \\
\hline
C5 & Code versioning system used & git \\
\hline
C6 & Software code languages, tools, and services used & Python, FEniCS, NumPy, PETSc, MPI, meshio, Gmsh \\
\hline
C7 & Compilation requirements, operating environments \& dependencies & FEniCS, meshio, Gmsh \\
\hline
C8 & If available Link to developer documentation/manual & \url{https://cashocs.readthedocs.io/} \\
\hline
C9 & Support email for questions & \href{mailto:sebastian.blauth@itwm.fraunhofer.de}{sebastian.blauth@itwm.fraunhofer.de} \\
\hline
\end{tabular}

\end{table}


\section{Description of the Software-Update}
\label{Description of the software-update}

This articles serves as update to our software cashocs \cite{Blauth2021cashocs} and discusses some of the key additions, changes, and improvements to the software. We describe both new functionalities which have been added to cashocs to extend its capabilities as well as structural improvements which make the software more efficient and easier to apply. 

\section{New Functionalities}
\label{sec:new_features}

Since its inception, cashocs has received many new functionalities, which we present in the following. For a detailed description of all changes and additions to cashocs, we refer to the release notes at \url{https://cashocs.readthedocs.io/en/stable/release_notes/}. One of the most practical additions to cashocs are space mapping methods which enable the solution of highly complex problems in the following way: The space mapping technique utilizes a model hierarchy consisting of a fine (detailed, complex) and a coarse (approximate, cheap) model. It provides an efficient way to optimize the fine model by successively optimizing and correcting a sequence of coarse model approximations. Particularly, there is no need for a direct optimization of the fine model. This makes the technique particularly interesting for industrial applications, where often commercial solvers without optimization capabilities are used for simulation. For more details on the space mapping technique, we refer the reader to, e.g., \cite{Bakr2001introduction,Echeverria2005Space,Blauth2023Space}. Moreover, we note that corresponding tutorials for using the space mapping technique in cashocs can be found at \url{https://cashocs.readthedocs.io/en/stable/user/demos/shape_optimization/demo_space_mapping_semilinear_transmission/} and \url{https://cashocs.readthedocs.io/en/stable/user/demos/shape_optimization/demo_space_mapping_uniform_flow_distribution/}.

\begin{figure}[!t]
	\centering
	\begin{subfigure}{\textwidth}
		\centering
		\begin{tikzpicture}
		\node at (0,0) {\includegraphics[width=0.55\textwidth, trim=0cm 0cm 0cm 46cm, clip]{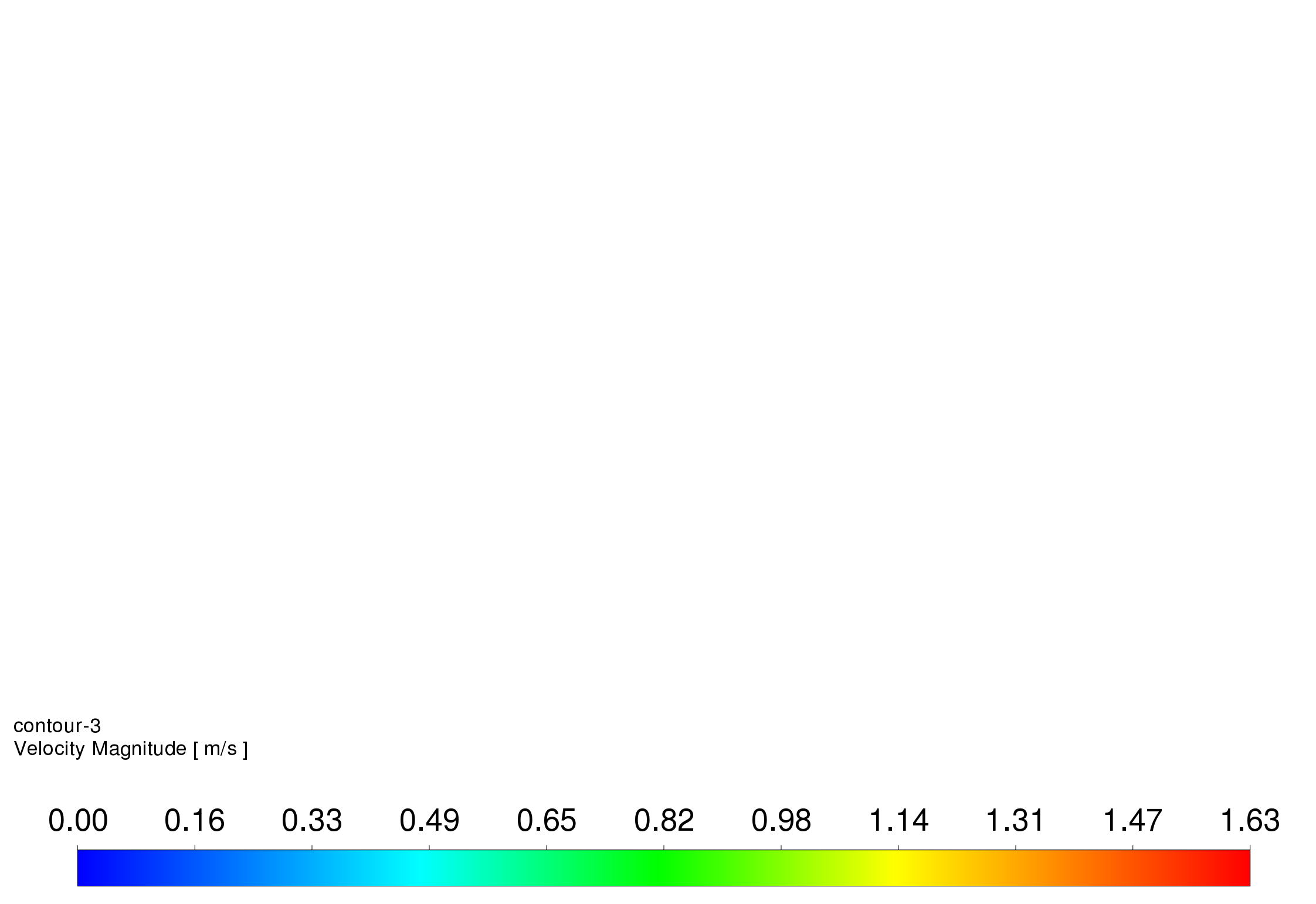}};
		
		\node at (0, 0.5) {\scriptsize Velocity Magnitude};
		\end{tikzpicture}
	\end{subfigure}

	\begin{subfigure}[!t]{0.5\textwidth}
		\centering
		\begin{minipage}[!t]{0.6\textwidth}
			\includegraphics[width=\textwidth]{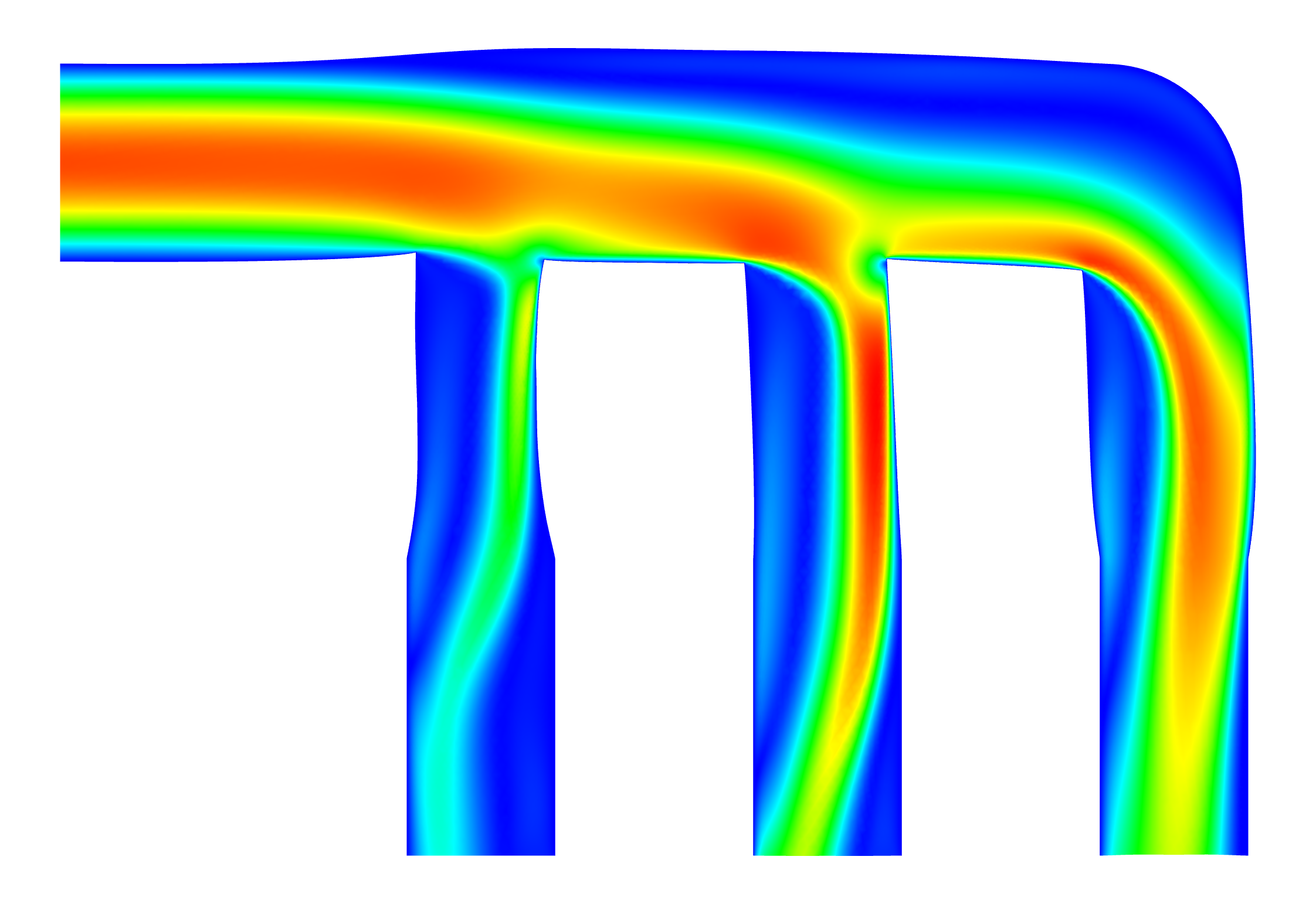}
		\end{minipage}%
		\hfil%
		\begin{minipage}[!t]{0.4\textwidth}
			\includegraphics[width=\textwidth]{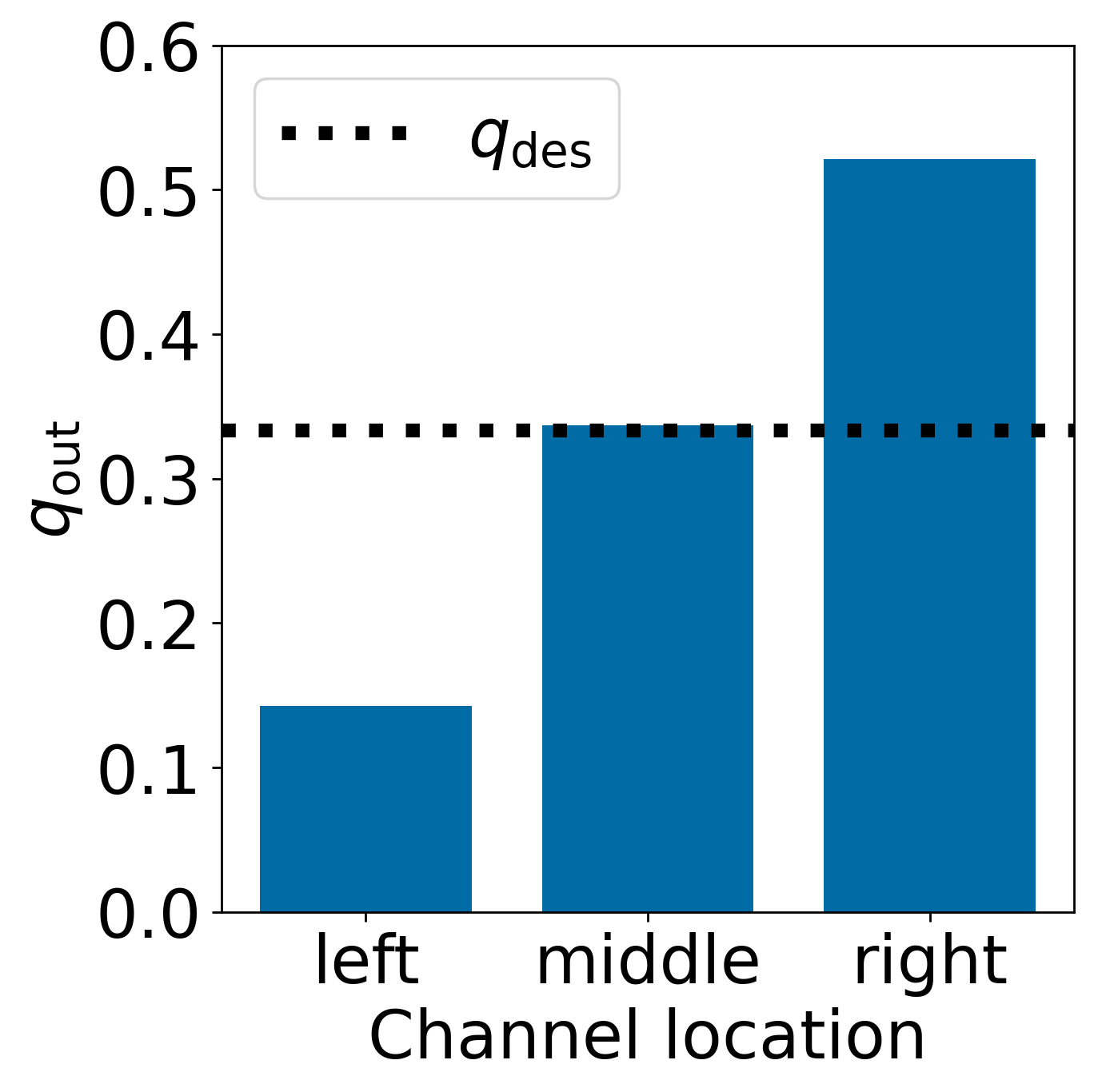}
		\end{minipage}
		\caption{Iteration 0.}
	\end{subfigure}%
	\hfil%
	\begin{subfigure}[!t]{0.5\textwidth}
		\centering
		\begin{minipage}[!t]{0.6\textwidth}
			\includegraphics[width=\textwidth]{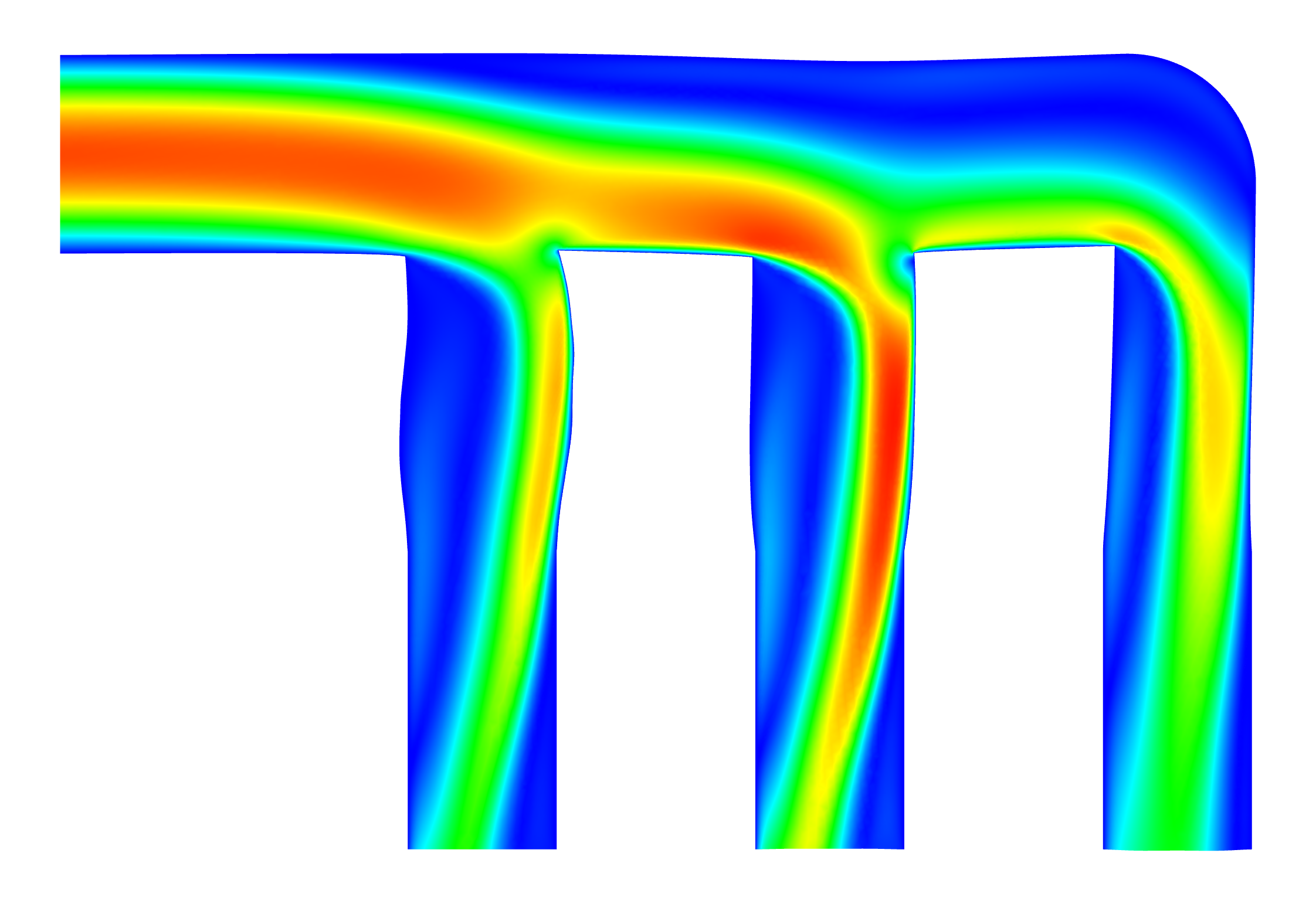}
		\end{minipage}%
		\hfil%
		\begin{minipage}[!t]{0.4\textwidth}
			\includegraphics[width=\textwidth]{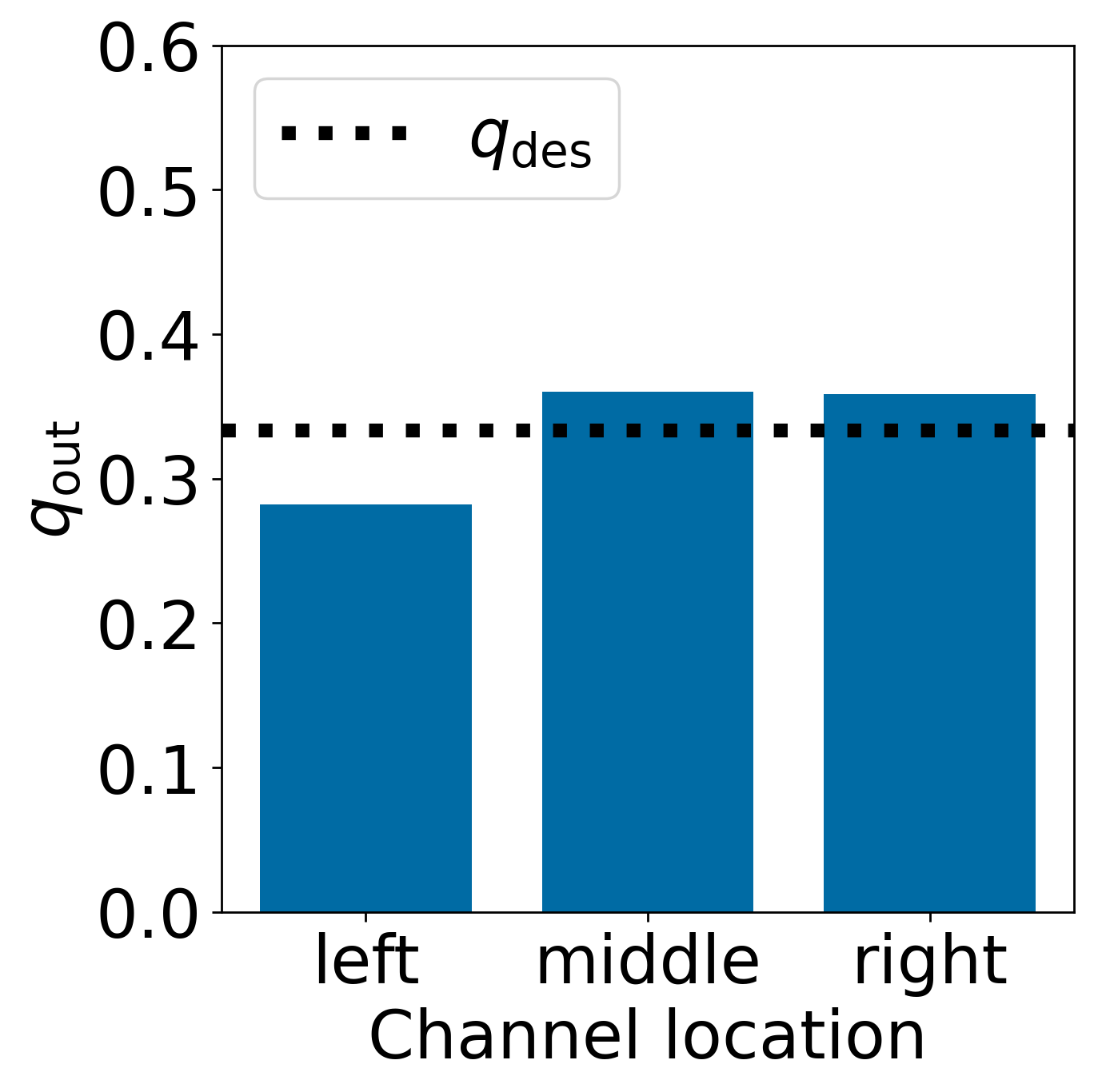}
		\end{minipage}
		\caption{Iteration 1.}
	\end{subfigure}

	\begin{subfigure}[!t]{0.49\textwidth}
		\centering
		\begin{minipage}[!t]{0.6\textwidth}
			\includegraphics[width=\textwidth]{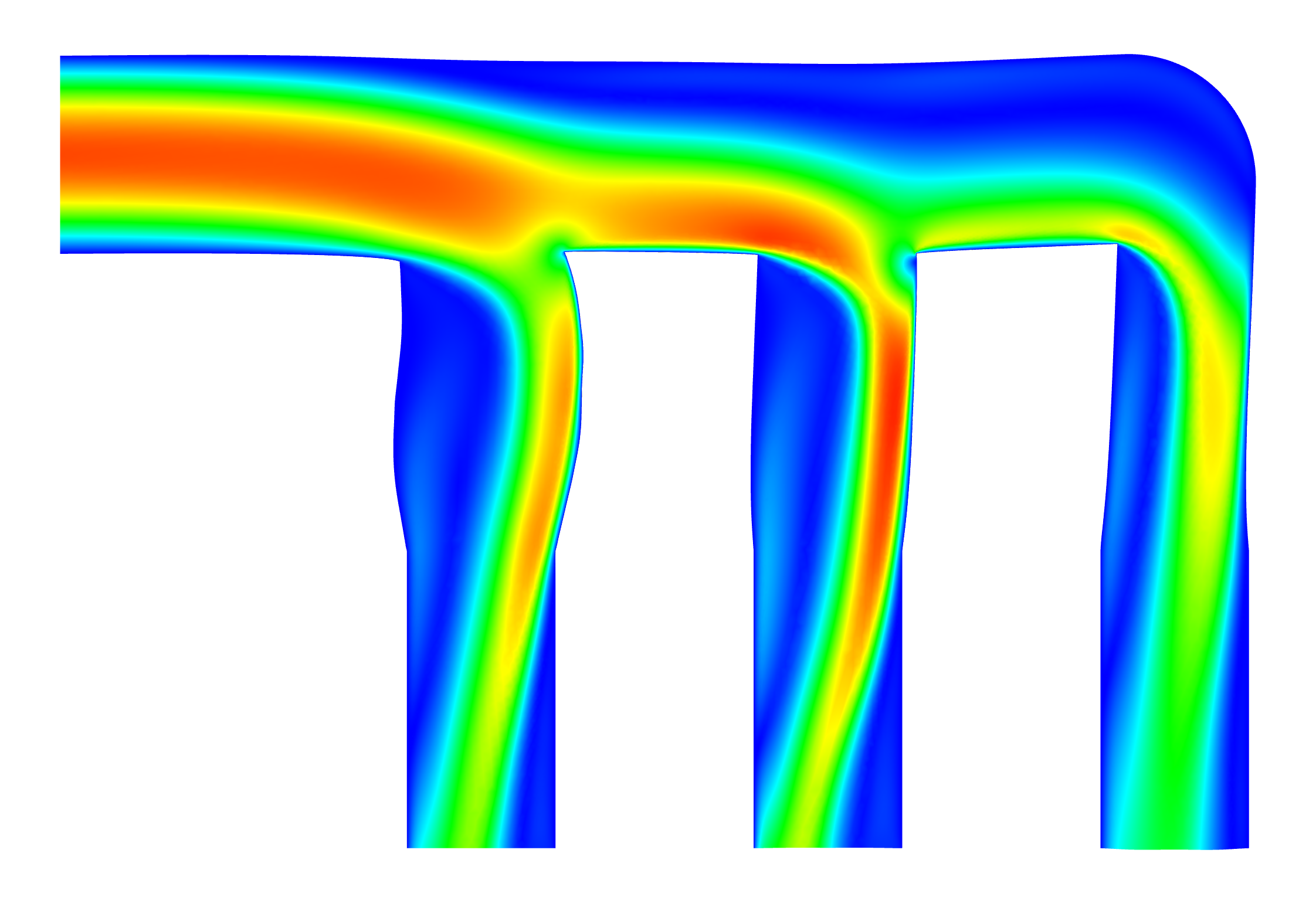}
		\end{minipage}%
		\hfil%
		\begin{minipage}[!t]{0.4\textwidth}
			\includegraphics[width=\textwidth]{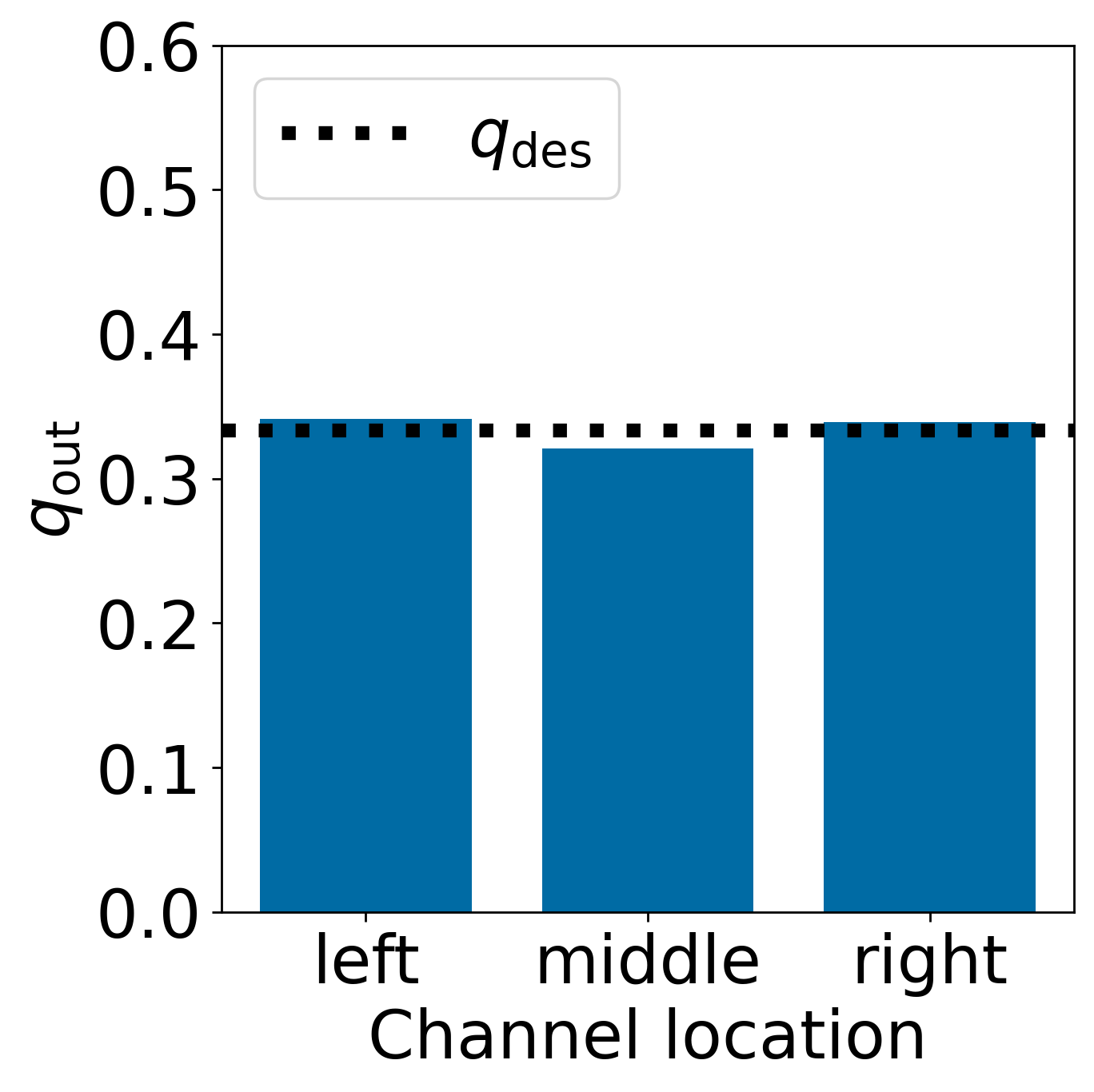}
		\end{minipage}
		\caption{Iteration 2.}
	\end{subfigure}%
	\hfil%
	\begin{subfigure}[!t]{0.49\textwidth}
		\centering
		\begin{minipage}[!t]{0.6\textwidth}
			\includegraphics[width=\textwidth]{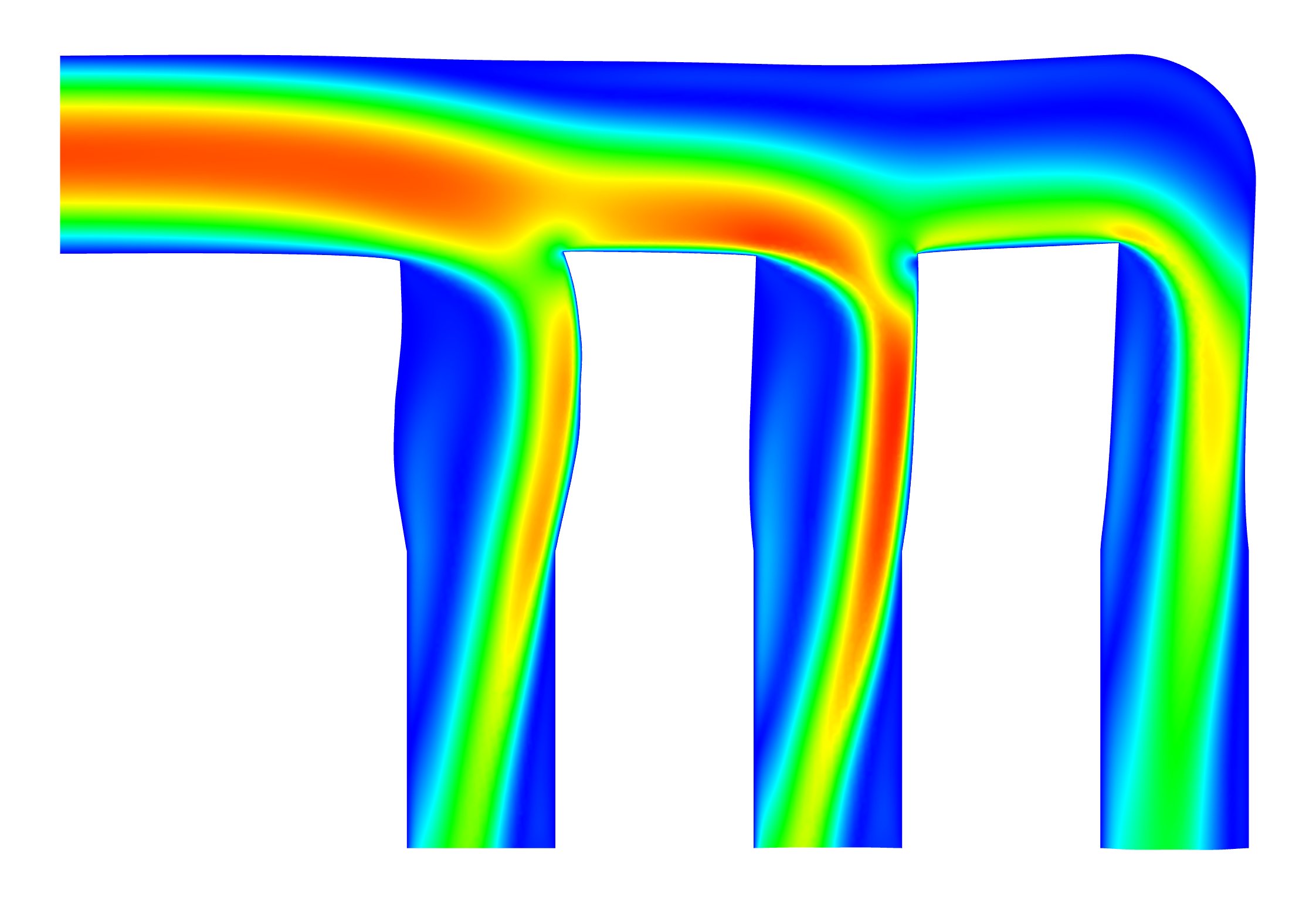}
		\end{minipage}%
		\hfil%
		\begin{minipage}[!t]{0.4\textwidth}
			\includegraphics[width=\textwidth]{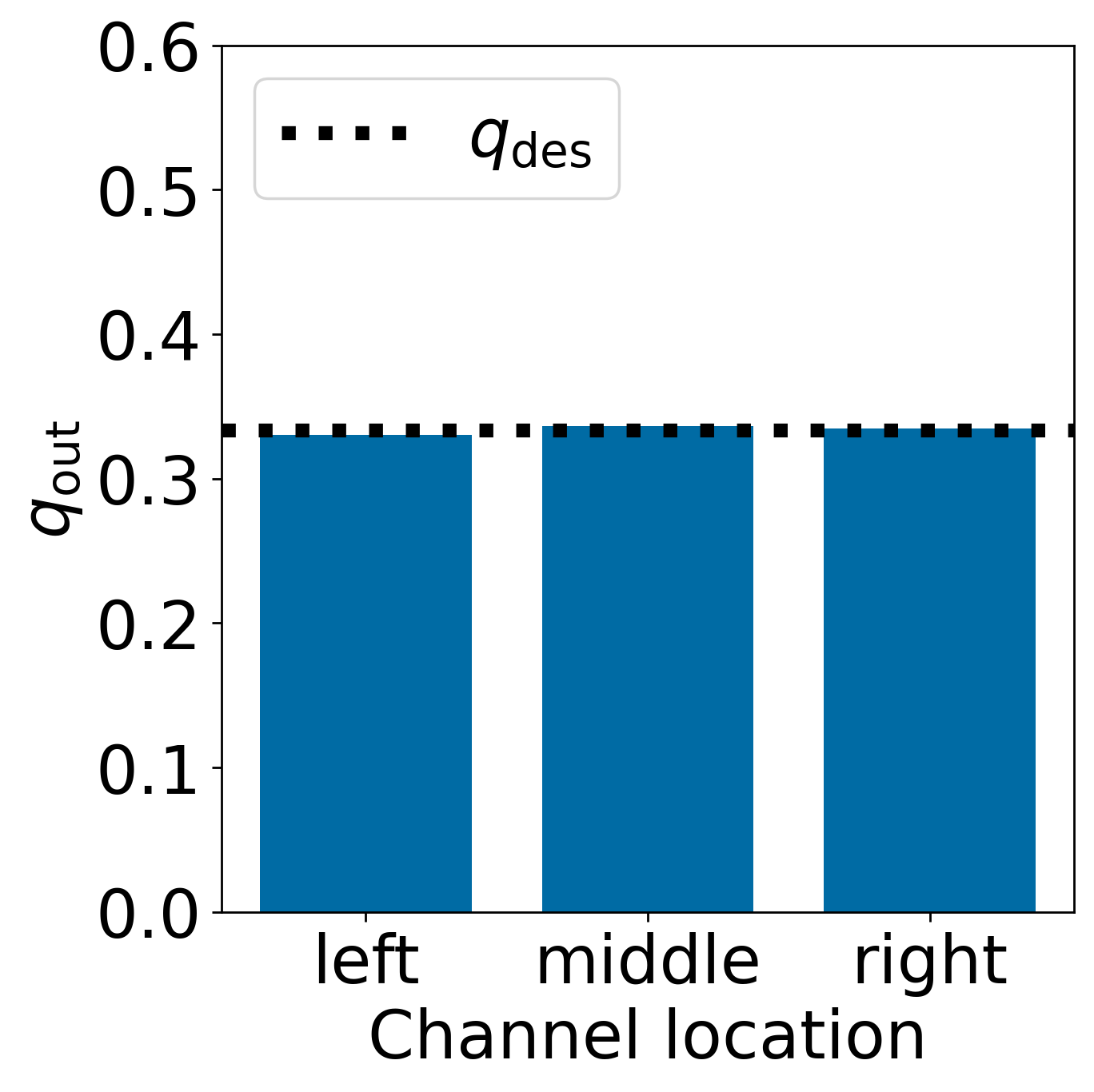}
		\end{minipage}
		\caption{Iteration 3.}
	\end{subfigure}

	\begin{subfigure}[!t]{0.49\textwidth}
		\centering
		\begin{minipage}[!t]{0.6\textwidth}
			\includegraphics[width=\textwidth]{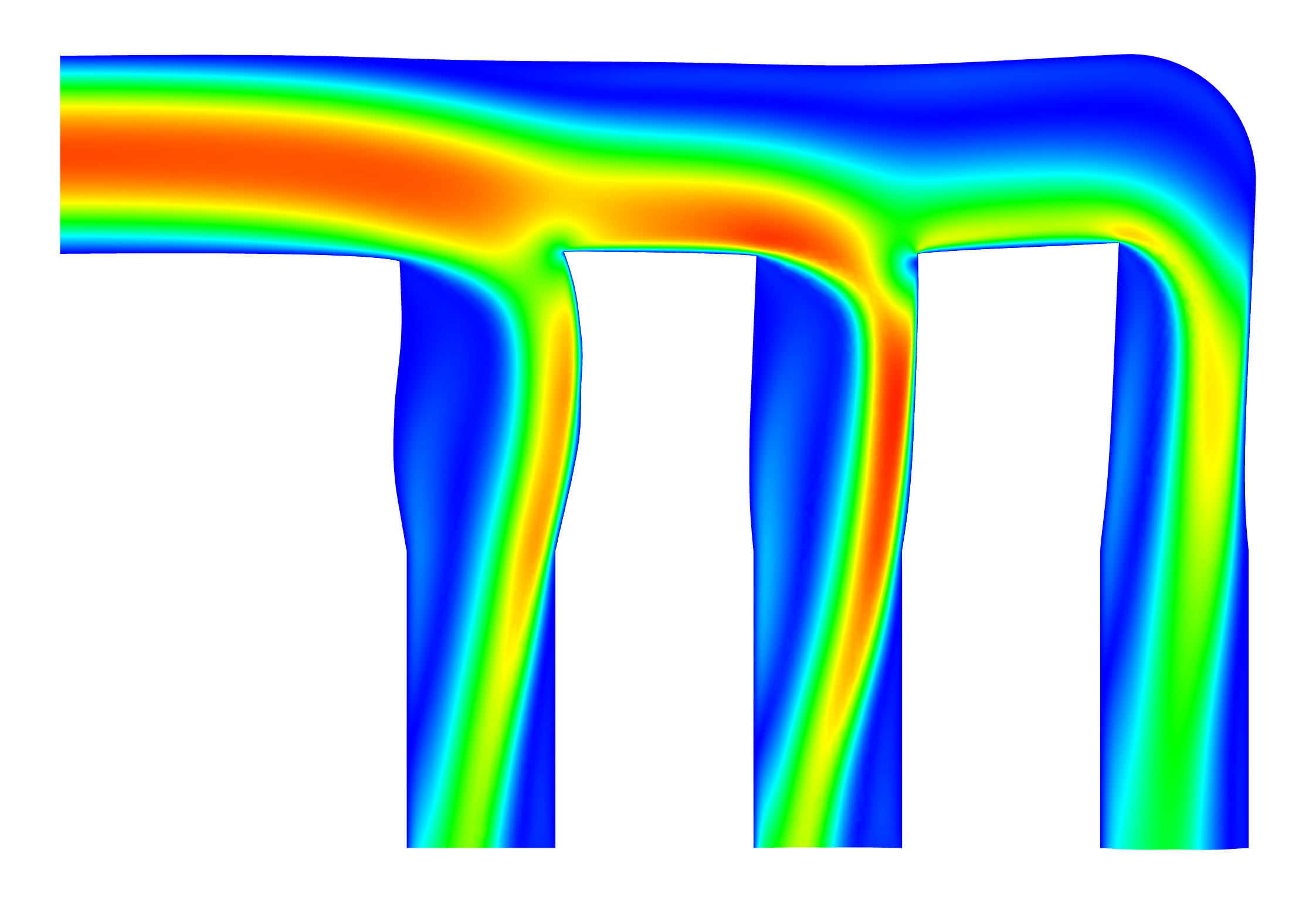}
		\end{minipage}%
		\hfil%
		\begin{minipage}[!t]{0.4\textwidth}
			\includegraphics[width=\textwidth]{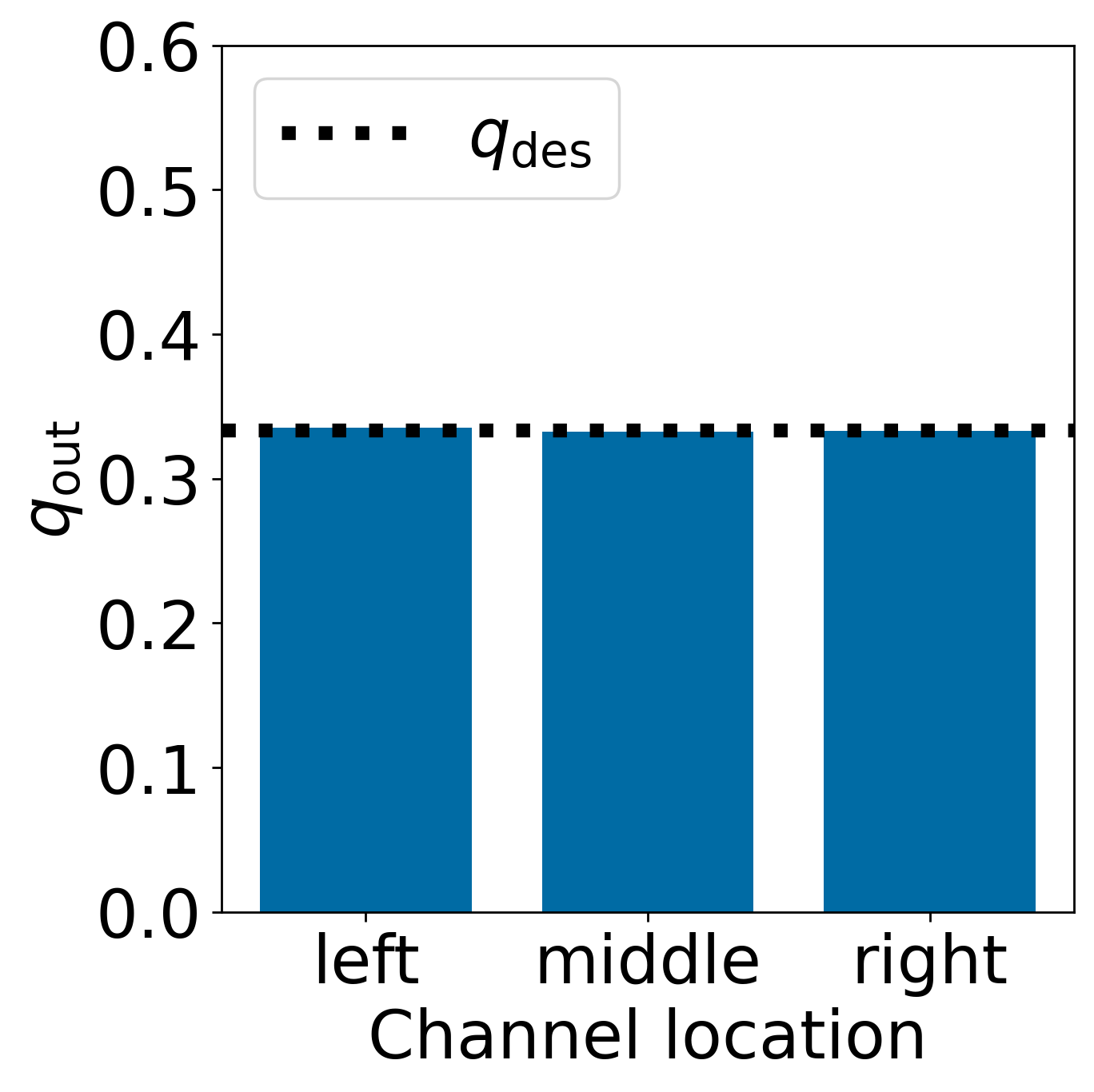}
		\end{minipage}
		\caption{Iteration 4.}
	\end{subfigure}
	\hfil
	\begin{subfigure}[!t]{0.49\textwidth}
		\centering
		\begin{minipage}[!t]{0.6\textwidth}
			\includegraphics[width=\textwidth]{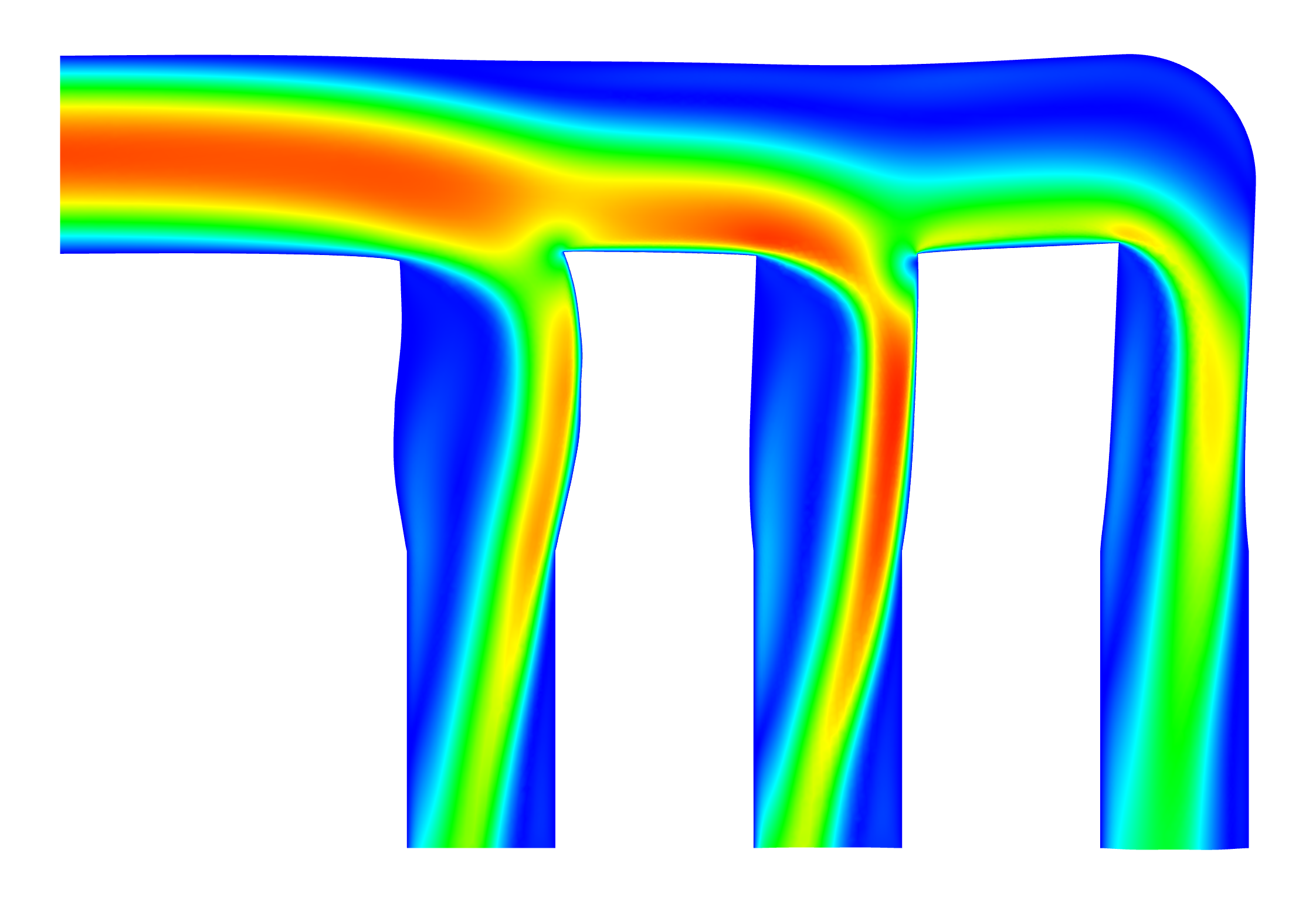}
		\end{minipage}%
		\hfil%
		\begin{minipage}[!t]{0.4\textwidth}
			\includegraphics[width=\textwidth]{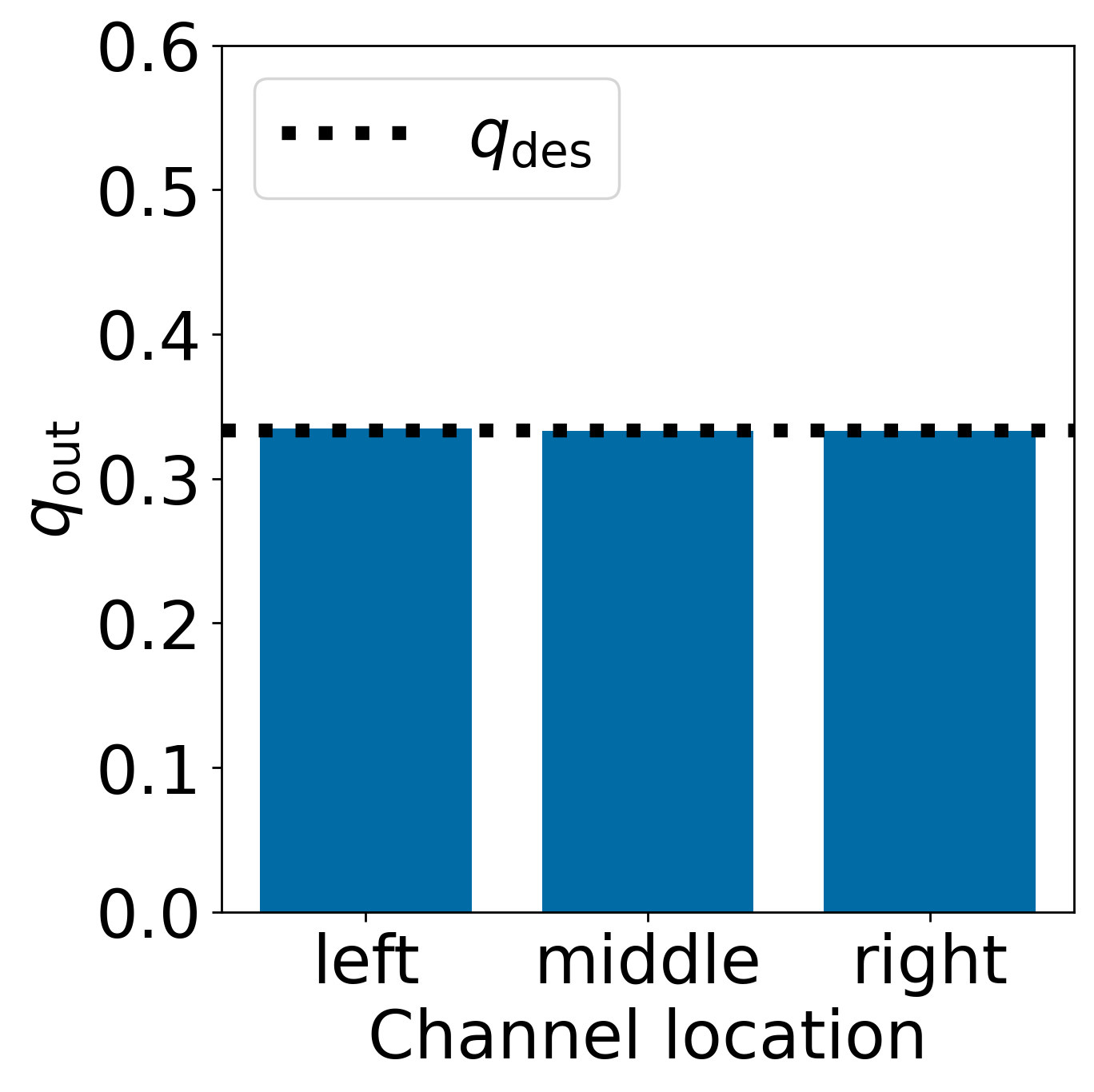}
		\end{minipage}
		\caption{Iteration 5.}
	\end{subfigure}
	
	\caption{Example for the space mapping method for shape optimization from \cite{Blauth2023Space}. The picture shows the evolution of the geometry, velocity field, and outlet flow rates over the space mapping iterations.}
	\label{fig:fluent_geometries}
\end{figure}

Whereas space mapping techniques for optimal control problems have already been investigated in the literature (see, e.g., \cite{Hintermueller2005Space,Marheineke2012Model,Totzeck2020Space}), space mapping methods for shape optimization have only been introduced very recently in \cite{Blauth2023Space}. Moreover, to the best of our knowledge, cashocs is the only software providing a framework for implementing and solving general space mapping problems in the context of PDE constrained optimization for both shape optimization and optimal control problems. 

As a reference, we briefly take a look at a problem from \cite{Blauth2023Space}, which is solved with the space mapping methods for shape optimization and is depicted in Figure~\ref{fig:fluent_geometries}. The aim of the optimization problem is to achieve a uniform flow distribution over three outlet pipes. For the fine model, the Navier-Stokes equations with a Reynolds number of 1000 are used and solved with Ansys Fluent \cite{Ansys2022Ansys}, whereas the coarse model is given by the linear Stokes system and is solved on a coarser grid with FEniCS \cite{Alnes2015FEniCS,Logg2012Automated,Logg2010DOLFIN}. We see that the space mapping method converges really fast in about five iterations and that the optimized geometry yields a nearly perfect uniform flow distribution. For more details, we refer the reader to \cite{Blauth2023Space} and to the corresponding tutorial at \url{https://cashocs.readthedocs.io/en/stable/user/demos/shape_optimization/demo_space_mapping_uniform_flow_distribution/}, where the corresponding tutorial can be found.

Another major new feature for cashocs is its support for solving topology optimization problems with a level-set method. Our software provides several optimization algorithms for such problems, including state-of-the-art algorithms for topology optimization \cite{Amstutz2006new} as well as novel quasi-Newton methods for topology optimization proposed in \cite{Blauth2023Quasi}. Due to the inherent difficulty of deriving topological derivatives, the latter cannot yet be computed with automatic differentiation methods but have to be supplied by the user. However, the automatic derivation of adjoint systems is already implemented in cashocs. This feature extends the applicability of cashocs to the field of topology optimization with topological sensitivity information. For more details on how to solve topology optimization problems with cashocs, we refer the reader to the corresponding tutorials at \url{https://cashocs.readthedocs.io/en/stable/user/demos/topology_optimization/}.

\begin{figure}[!t]
	\centering
	\tikzstyle{block} = [rectangle, draw, text width=12em, text centered, rounded corners, minimum height=4em]
	\tikzstyle{splitblock} = [rectangle split, draw, text centered, rounded corners, minimum height=4em, rectangle split parts=2]
	\tikzstyle{line} = [draw, -latex']
	
	\begin{tikzpicture}[scale=0.2]
	\node[splitblock](vanilla) at (0, 0) {\textbf{PDE Constrained Optimization Problems} \nodepart[align=left]{two} \texttt{cashocs.ShapeOptimizationProblem} \\ \texttt{cashocs.OptimalControlProblem} \\ \texttt{\color{blue} cashocs.TopologyOptimizationProblem} };
	
	\node[splitblock, below of=vanilla, node distance=6em](constrained) {\textbf{\color{blue} Additional Constraints} \nodepart[align=left]{two} \texttt{\color{blue} cashocs.ConstrainedShapeOptimizationProblem} \\ \texttt{\color{blue} cashocs.ConstrainedOptimalControlProblem} };
	
	\node[splitblock, below of=constrained, node distance=6em](space_mapping) {\textbf{\color{blue} Space Mapping Problems} \nodepart[align=left]{two} \texttt{\color{blue} cashocs.space\_mapping.shape\_optimization.SpaceMappingProblem} \\ \texttt{\color{blue} cashocs.space\_mapping.optimal\_control.SpaceMappingProblem} };
	
	\end{tikzpicture}
	\caption{Overview of the types of optimization problems that can be solved with cashocs. Items colored in blue are new additions to version 2.0.}
	\label{fig:architecture}
\end{figure}
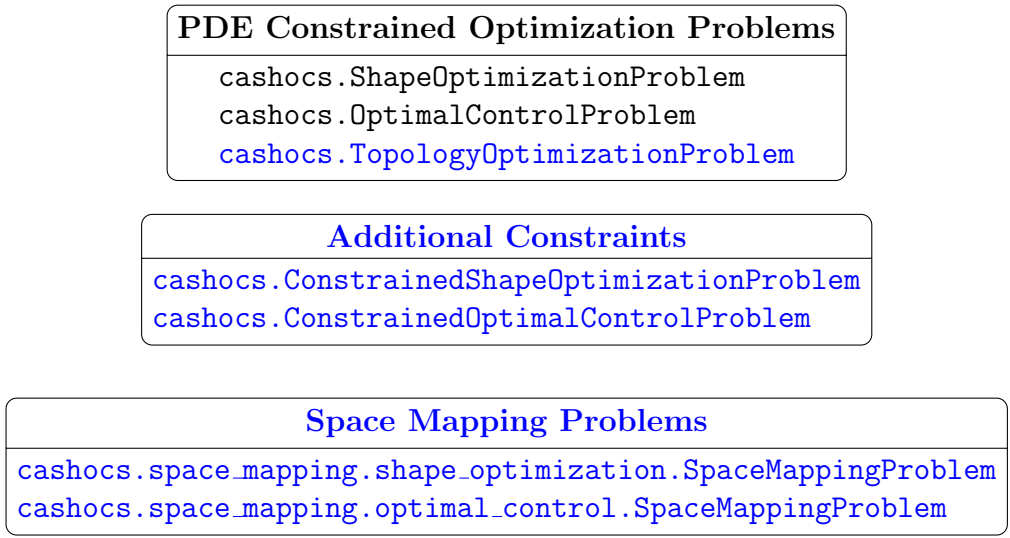

Further, cashocs has received methods for the automatic treatment of additional constraints for optimization problems, such as state or control constraints. In particular we have implemented a quadratic penalty and an augmented Lagrangian method (see, e.g., \cite{Nocedal2006Numerical}) for dealing with constraints other than PDEs. These methods are applicable for a wide variety of constraints, including state constraints, which opens up new problem classes for cashocs to solve. We refer the reader to \url{https://cashocs.readthedocs.io/en/stable/user/demos/optimal_control/demo_constraints/}, where a demo for using additional constraints with cashocs can be found. An overview over all types of optimization problems that can be solved with cashocs is given in Figure~\ref{fig:architecture}.

The new version of cashocs also comes with the ability to automatically scale individual terms of the cost functional based on their magnitude for the initial iteration. This makes it easier for users to weigh different parts of the cost functional according to their needs. Moreover, this also facilitates the solution of multi-criteria optimization problems with scalarization methods, see, e.g., \cite{Ehrgott2005Multicriteria,Blauth2023Multi}. We refer the reader to the corresponding cashocs tutorial at \url{https://cashocs.readthedocs.io/en/stable/user/demos/shape_optimization/demo_scaling/}.

With the new version of cashocs, users now have the possibility to define custom scalar products for the optimization, which is particularly useful for shape optimization, where the choice of the scalar product for computing the shape gradient is particularly important. The reader is referred to the cashocs tutorials at \url{https://cashocs.readthedocs.io/en/stable/user/demos/shape_optimization/demo_scaling/} for more information. Moreover, the updated version of cashocs also supports computing the shape gradient with the $p$-Laplace equations, based on the approach presented in \cite{Mueller2021novel}. This approach may preserve the mesh quality and can yield better optimized geometries for certain problems. For more details, the reader is referred to the corresponding demo at \url{https://cashocs.readthedocs.io/en/stable/user/demos/shape_optimization/demo_p_laplacian/}.

Finally, we have also added more sophisticated line search methods to cashocs, which make use of polynomial models to efficiently compute a suitable stepsize. During the line search, a quadratic or cubic model of the cost functional along the search direction is formed and minimized to obtain suitable stepsizes in fewer iterations compared to a classical Armijo backtracking search (see, e.g., \cite{Kelley1999Iterative}).

\section{Structural Improvements}
\label{sec:technical_features}

In addition to the new features described above, there have also been some major structural improvements to cashocs which we discuss in the following.

The most important structural improvement is the support for parallelism with MPI. In particular, most of the PDE constrained optimization problems that can be solved in serial can now also be solved in parallel without the need for any code modifications. Therefore, cashocs is now able to be run on high-performance-computing systems and, thus, can solve huge optimization problems which cannot be treated in serial. The MPI implementation builds on and leverages the MPI support of FEniCS, so that there are only minimal code changes necessary on the user's side. 

Another structural improvement for cashocs is a new remeshing workflow which overcomes the limitations induced by the previous one. The benefit of the new workflow is that remeshing can now be applied also for the space mapping and constraint handling problems, which was not possible beforehand. As remeshing is of particular importance for solving shape optimization problems, this change greatly increases the applicability of cashocs. However, user's must now follow a slightly different, but still straightforward, syntax to solve shape optimization problems with remeshing. The corresponding tutorial can be found at \url{https://cashocs.readthedocs.io/en/stable/user/demos/shape_optimization/demo_remeshing/}. 

Finally, we note that the remeshing workflow itself is not parallelized with MPI due to limitations of the used software tools Gmsh, which supports only OpenMP parallelism, and meshio, which does not support parallelism at all. This can potentially create bottlenecks when solving problems where lots of remeshing is required. However, as usually the solution of the state and adjoint systems is the most CPU and memory consuming task, the bottleneck of the remeshing is not too severe for many applications.

\section{Conclusion and Outlook}

With the update to version 2.0, our software cashocs has gained significant new features, such as a space mapping framework, support for parallelism via MPI, and the ability to treat topology optimization problems. These additions make cashocs even more flexible and relevant for solving PDE constrained optimization problems for practical and industrial applications. 



\bibliographystyle{elsarticle-num} 
\bibliography{literature_db.bib}

\end{document}